\documentclass[smallextended,numbook,runningheads]{svjour3}     % onecolumn (second format)
\smartqed  % flush right qed marks, e.g. at end of proof
\usepackage{graphicx}
\usepackage{epsfig}
\usepackage{amsmath}
\usepackage{mathptmx}      % use Times fonts if available on your TeX system
\journalname{...}
\begin{document}

\title{High phase-lag order trigonometrically fitted two-step Obrechkoff methods
for the numerical solution of periodic initial value problems%\thanks{Grants or other notes about the
%article that should go on the front page should be placed here.
% General acknowledgments should be placed at the end of the article.
}%}
%\subtitle{Do you have a subtitle\\ If so, write it here}

\titlerunning{High phase-lag order trigonometrically fitted two-step Obrechkoff methods}        % if too long for running head

\author{Ali Shokri\and
        Hosein Saadat
}


\institute{Ali Shokri \at
              Faculty of Mathematical Science, University of
              Maragheh, Maragheh, Iran. \\
              \email{shokri@maragheh.ac.ir}           %  \\
           \and
           Hosein Saadat \at
              Faculty of Mathematical Science, University of
              Maragheh, Maragheh, Iran.\\
              \email{hosein67saadat@yahoo.com}
}

\date{Received: date / Accepted: date}

\maketitle

\begin{abstract}
In this paper, we present the two-step trigonometrically fitted
symmetric Obrechkoff methods with algebraic order of twelve. The
method is based on the symmetric two-step Obrechkoff method, with 12
algebraic order, high phase-lag order and is constructed to solve
IVPs with periodic solutions such as orbital problems. We compare
the new method to some recently constructed optimized methods from
the literature. The numerical results obtained by the new method for
some problems show its superiority in efficiency, accuracy and
stability.
%Include keywords and mathematical subject classification numbers as needed.
\subclass{MSC 65l05 \and MSC 65l07 \and 65l20}
\end{abstract}
%\\\\\\\\\\\\\\\\\\\\\\\\\\\\\\\\\\\\\\\\\\\\\\\\\\\\\\\\\\\\\\\\\\\\\\\\\\\\\\\\\\\\\\\\\\\\\\\\\\\\\\\\\\\\
\section{Introduction}
\label{intro}
In this paper, the symmetric Obrechkoff methods for
solving special class of initial value problems associated with
second order ordinary differential equations of the type
\begin{equation} \label{dif}
y''=f(x,y),\quad y(x_0)=y_0,\quad y'(x_0)=y'_0,
\end{equation}
in which the first order derivatives do not occur explicitly, are
discussed.The numerical integration methods for (\ref{dif}) can be
divided into two distinct classes:
\begin{enumerate}
\item Problems for which the solution period is known (even approximately) in advance.
\item Problems for which the period is not known.
\end{enumerate}
For several decades, there has been strong interest in searching for
better numerical methods to integrate first-order and second-order
initial value problems, because these problems are usually
encountered in celestial mechanics, quantum mechanical scattering
theory, theoretical physics and chemistry, and electronics.
Generally, the solution of $(1)$ is periodic, so it is expected that
the result produced by some numerical methods preserves the
analogical periodicity of the analytic solution [9-22].
Computational methods involving a parameter proposed by Gautschi
\cite{Ga}, Jain et al. \cite{Ja}, Sommeijer and et al \cite{So} and Steifel and Bettis \cite{St}
yield numerical solution of problems of class (1). Chawla and et al
\cite{cha,ch,chaw}, Ananthakrishnaiah \cite{Ana}, Shokri and et al.
\cite{Sho1,Sh,Sho2,Sho}, Dahlquist \cite{Dal}, Franco \cite{Fra}, Lambert and
Watson [9], Tsitouras and Simos \cite{Tsi}, Simos and et al. \cite{Si,Sim,S}, Hairer \cite{Har}, Wang
et al. \cite{W,Wa,Wan}, Saldanha and Achar \cite{Sa}, and Daele
and Vanden Berghe \cite{Da} have developed methods to solve problems
of class (2). Consider Obrechkoff method of the form
\begin{equation} \label{ob}
\sum _{i=0}^{k}\alpha _iy_{n-j+1}=\sum _{i=1}^lh^{2i}\sum _{j=0}^k
\beta_{ij} y_{n-j+1}^{(2i)},
\end{equation}
for the numerical integration of the problem (\ref{dif}). The method
(\ref{ob}) is symmetric when $\alpha_j = \alpha _{k-j}$ , $\beta_j =
\beta_{k-j}$ , $ j=0,1,2,\cdots , k$ , and it is of order $ q $ if
the truncation error associated with the linear difference operator
is given as
\[TE=C_{q+2}h^{q+2}y^{(q+2)},\qquad x_{n-k+1}<\eta< x_{n+1},\]
where $C_{q+2}$ is a constant dependent on $h$. When the method
(\ref{ob}) is applied to the test problem, we get the characteristic
equation as
\begin{equation}\label{ro}
\rho(\xi)-\sum_{i=1}^l(-1)^iv^{2i}\sigma_i(\xi)=0,
\end{equation}
where $v=\lambda h$ and
\begin{equation} \label{6}
\rho(\xi)=\sum_{j=0}^k\alpha_j\xi^{k-j},\quad\sigma_i(\xi)=
\sum_{j=0}^k\beta_{ij}\xi^{k-j},\quad i=1,2,\cdots,l.
\end{equation}
\begin{definition}
The method (\ref{ob}) is said to have interval of periodicity
$(0,v_0^2)$ if for all $v^2\in(0,v_0^2)$ the roots of Eq. (\ref{ro})
are complex and at least two of them lie on the unit circle and the
others lie inside the unit circle.
\end{definition}
\begin{definition}
The method (\ref{ob}) is said to be P-stable if its interval of
periodicity is $(0,\infty)$.
\end{definition}
\begin{definition}
For any symmetric multistep methods, the phase-lag (frequency
distortion) of order $q$ is given by
\begin{equation} \label{7}
t(v)=v-\theta(v)=Cv^{q+1}+O(v^{q+2}),
\end{equation}
where $C$ is the phase lag constant and $q$ is the phase-lag order.
\end{definition}
The characteristic equation of the method (\ref{ob}) is given by
\begin{equation}\label{char}
\Omega(s:v^2)=A(v)s^2-2B(v)s+A(v)=0,
\end{equation}
where
\begin{equation} \label{14}
A(v)=1+\sum_{i=1}^m(-1)^i\beta_{i0}v^{2i},\quad
B(v)=1+\sum_{i=1}^m(-1)^i\beta_{i1}v^{2i},
\end{equation}
$\Psi$ contains polynomial functions together with trigonometric
polynomials
\begin{equation}
\Psi_{trig}=\left\{1,t,\cdots,t^K,\cos(r\omega t),\sin(r\omega
t),\quad\quad r=1,2,\cdots,P\right\}.
\end{equation}
The resulting methods are then based on a hybrid set of polynomials
and trigonometric functions. If $P$ is limited to $P=\frac{M-2}{2}$,
we called method with zero phase-lag.
\begin{remark}
We present here the trigonometric versions of the set. In case
$\omega$ is purely imaginary one obtains the hyperbolic description
of this set. This set is characterized by two integer parameters
$K$ and $P$. The set in which there is no polynomial part is
identified by $K=-1$ while the set in which there is no
trigonometric polynomial component is identified by $K=-1$. For each
problem one has $K+2P=M-3$, where $M-1$ is the maximum exponent
present in the full polynomial basis for the same problem.
\end{remark}

%\\\\\\\\\\\\\\\\\\\\\\\\\\\\\\\\\\\\\\\\\\\\\\\\\\\\\\\\\\\\\\\\\\\\\\\\\\\\\\\\\\\\\\\\\\\\\\\\\\\\\\\\\\\\\\\\
\section{Construction of the new method}
\label{sec:1} From the form \eqref{ob} and without loss of
generality we assume $\alpha_j=\alpha_{m-j}$,
$\beta_{i,j}=\beta_{i,m-j}$, $j=0(1)\lfloor \frac{m}{2}\rfloor$ and
we can write
\begin{equation} \label{ob2}
y_{n+1}-2y_{n}+ y_{n-1}=\sum_{i=1}^mh^{2i}\left[\beta_{i0}
y_{n+1}^{(2i)}+\beta_{i1}y_{n}^{(2i)}+\beta _{i 0}
y_{n-1}^{(2i)}\right],
\end{equation}
when $m=3$ we get
\begin{eqnarray}\label{ob3}
y_{n+1}-2y_{n}+y_{n-1}&=&h^{2}\left[\beta _{10}
(y_{n+1}^{(2)}+y_{n-1}^{(2)})+\beta _{11}
y_{n}^{(2)}\right]\nonumber\\
&+&h^4\left[\beta_{20}(y_{n+1}^{(4)}+y_{n-1}^{(4)})+\beta
_{21}y_n^{(4)}\right]\nonumber\\
&+&h^6\left[\beta_{30}(y_{n+1}^{(6)}+y_{n-1}^{(6)})+\beta
_{31}y_n^{(6)}\right].
\end{eqnarray}
$M-3$ for method (\ref{ob3}) is 11 so that if $P=-1$, $K=13$ we
obtain classic method and the coefficients of this method are
\begin{eqnarray}\label{cla}
\beta_{1,0}&=&\frac{229}{7788},\quad\beta_{1,1}=\frac{3665}{3894},\quad\beta_{2,0}=-\frac{1}{2360},\nonumber\\
\beta_{2,1}&=&\frac{711}{12980},\quad\beta_{3,0}=\frac{127}{39251520},\quad\beta_{3,1}=\frac{2923}{3925152},
\end{eqnarray}
where its phase-lag is given by
\[pl_{clas}:=-\frac{45469}{3394722659328000}v^{12}+O\left(v^{14}\right),\]
and its local truncation error is given by
\[LTE_{clas}=-\frac{45469}{1697361329664000}h^{14}y^{(14)}+O\left(h^{16}\right).\]
If $P=6$, $ K=-1$ we obtain the method with zero phase-lag (PL), and
the coefficients of this case are given in \cite{Sa}.
%%%%%%%%%%%%%%%%%%%%%%%%%%%%%%%%%%%%%%%%%%%%%%%%%%%%%%%%%%%%%%%%%%%%%%%%%%%%%%%
\subsection{The first formula}
If $P=0$, $K=11$, so we called PL$'$, we have
\[\beta_{1,0}=\frac{1}{6v^2}\frac{\beta_{1,0num}}{A},\quad
\beta_{1,1}=\frac{1}{3v^2}\frac{\beta_{1,1num}}{A},\quad
\beta_{2,0}=\frac{-1}{5040v^2}\frac{\beta_{2,0num}}{A},\]
\begin{equation}\label{co3}
\beta_{2,1}=\frac{1}{2520v^2}\frac{\beta_{2,1num}}{A},\quad
\beta_{3,0}=\frac{-1}{10080v^2}\frac{\beta_{3,0num}}{A},\quad
\beta_{3,1}=\frac{1}{5040v^2}\frac{\beta_{3,1num}}{A},
\end{equation}
where
\[A=15120\cos v-15120+6900v^2-313v^4+660v^2\cos v+13v^4\cos v,\]
and
\begin{eqnarray}
\beta_{1,0num}&=&-45360v^2+3702v^4-89v^6+78v^4\cos v+2v^6\cos
v+90720-90720\cos v,\nonumber\\
\beta_{1,1num}&=&45360v^2\cos v+16998v^4-850v^6+37v^6\cos
v-90720+90720\cos v\nonumber\\
&+&1902v^4\cos v,\nonumber\\
\beta_{2,0num}&=&-65520v^2\cos
v-1597680v^2+105840v^4-1907v^6+17v^6\cos v+3326400\nonumber\\
&-&3326400\cos
v,\nonumber\\
\beta_{2,1num}&=&3109680v^2\cos v+14278320v^2-30257v^6+1907v^6\cos
v\nonumber\\
&-&34776000+34776000\cos v+105840v^4\cos v,\nonumber\\
\beta_{3,0num}&=&3360v^2\cos v+62160v^2-3814v^4+59v^6+34v^4\cos
v-131040+131040\cos v,\nonumber\\
\beta_{3,1num}&=&149520v^2\cos v+1428000v^2-60514v^4+59v^6\cos
v-3155040\nonumber\\
&+&3155040\cos v+3814v^4\cos v,\nonumber
\end{eqnarray}
for small values of $v$ the above formulae are subject to heavy
cancelations. In this case the following Taylor series expansion
must be used:
\begin{eqnarray}
\beta_{1,0}&=&\frac{229}{7788}+\frac{45469}{1314147120}v^2+\frac{85771}{341152592352}v^4+
\frac{42739761203}{29358705101073004800}v^6\nonumber\\
&+&\frac{3801508031029}{608197283570236453277184}v^8+\frac{168279971604233}{13575027728788584540475136000}v^{10}\nonumber\\
&-&\frac{266348222900207221}{2703381808485285252094734713548800}v^{12}+\cdots,\nonumber
\end{eqnarray}
\begin{eqnarray}
\beta_{1,1}&=&\frac{3665}{3894}-\frac{45469}{657073560}v^2-\frac{85771}{170576296176}v^4-
\frac{42739761203}{14679352550536502400}v^6\nonumber\\
&-&\frac{3801508031029}{304098641785118226638592}v^8
-\frac{168279971604233}{6787513864394292270237568000}v^{10}\nonumber\\
&+&\frac{266348222900207221}{1351690904242642626047367356774400}v^{12}+\cdots,\nonumber
\end{eqnarray}
\begin{eqnarray}
\beta_{2,0}&=&-\frac{1}{2360}-\frac{45469}{30105915840}v^2-\frac{12253}{1116499393152}v^4-
\frac{42739761203}{672581244133672473600}v^6\nonumber\\
&-&\frac{3801508031029}{13933246859972689656895488}v^8
-\frac{168279971604233}{310991544332247573109066752000}v^{10}\nonumber\\
&+&\frac{266348222900207221}{61932019612571989411624831619481600}v^{12}+\cdots,\nonumber
\end{eqnarray}
\begin{eqnarray}
\beta_{2,1}&=&\frac{711}{12980}-\frac{1045787}{33116507424}v^2-\frac{1409095}{6140746662336}v^4
-\frac{983014507669}{739839368547039720960}v^6\nonumber\\
&-&\frac{437173423568335}{76632857729849793112925184}v^8
-\frac{3870439346897359}{342090698765472330419973427200}v^{10}\nonumber\\
&+&\frac{266348222900207221}{2961966155383877754469013686149120}v^{12}+\cdots,\nonumber
\end{eqnarray}
\begin{eqnarray}
\beta_{3,0}&=&\frac{127}{39251520}+\frac{45469}{1528454188800}v^2+\frac{12253}{56683815344640}v^4
+\frac{42739761203}{34146432394478756352000}v^6\nonumber\\
&+&\frac{3801508031029}{707380225198613474888540160}v^8+
\frac{168279971604233}{15788801481483338327075696640000}v^{10}\nonumber\\
&-&\frac{266348222900207221}{3144240995715193308590183759142912000}v^{12}+\cdots,\nonumber
\end{eqnarray}
\begin{eqnarray}
\beta_{3,1}&=&\frac{2923}{3925152}-\frac{14231797}{9934952227200}v^2
-\frac{3835189}{368444799740160}v^4-\frac{13377545256539}{221951810564111916288000}v^6\nonumber\\
&-&\frac{1189872013712077}{4597971463790987586775511040}v^8-
\frac{52671631112124929}{102627209629641699125992028160000}v^{10}\nonumber\\
&+&\frac{83366993767764860173}{20437566472148756505836194434428928000}v^{12}+\cdots.\nonumber
\end{eqnarray}
The phase-lag and the local truncation error for the PL$'$ method
are given by
\begin{eqnarray}
LTE_{PL'}&=&(1-\beta_{1,1}-2\beta_{1,0})h^2y_n^{(2)}+
\left(\frac{1}{12}-\beta_{1,0}-2\beta_{2,0}-\beta_{2,1}\right)h^4y^{(4)}\nonumber\\
&+&\left(\frac{1}{360}-\frac{\beta_{1,0}}{12}-\beta_{2,0}-2\beta_{3,0}-\beta_{3,1}\right)h^6y_n^{(6)}
+\left(\frac{2}{8!}-\frac{2\beta_{1,0}}{6!}-\frac{2\beta_{2,0}}{4!}-\frac{2\beta_{3,0}}{2!}\right)h^8y_n^{(8)}\nonumber\\
&-&\left(\frac{2}{10!}-\frac{2\beta_{1,0}}{8!}-\frac{2\beta_{2,0}}{6!}-\frac{2\beta_{3,0}}{4!}\right)h^{10}y^{(10)}
+\left(\frac{2}{12!}-\frac{2\beta_{1,0}}{10!}-\frac{2\beta_{2,0}}{8!}-\frac{2\beta_{3,0}}{6!}\right)h^{12}y^{(12)}\nonumber\\
&+&\left(\frac{2}{14!}-\frac{2\beta_{1,0}}{12!}-\frac{2\beta_{2,0}}{10!}-\frac{2\beta_{3,0}}{8!}\right)h^{14}y^{(14)}
+O\left(h^{16}\right).\nonumber
\end{eqnarray}
hence
\[pl_{PL'}=\frac{731602960042513638469539403}{1287287007659726361217210431335975522416459776000000}v^{24},\]
and
\[LTE_{PL'}=-\frac{45469}{1697361329664000}\left(y^{(14)}+\omega^2y^{(12)}\right)h^{14},\]
where $v=\omega h$, $\omega$ is the frequency and $h$ is the step
length. As $v\rightarrow0$, the LTE of the method \eqref{ob3} with
derived coefficients \eqref{co3} tends to
$\frac{45469}{169736132966400}h^{14}y^{(14)}+O\left(h^{16}\right)$.
which agrees with the LTE of the three methods due to Wang
\cite{Wan}, Simos \cite{Si} and Daele \cite{Da}, Achar \cite{Achar},
as $H\rightarrow0$. The behavior of the coefficients of the PL$'$
method are shown in Figures 2.1, to 2.6.
\begin{figure} \epsfxsize=7cm\epsfysize=7cm
\begin{center}\epsfbox{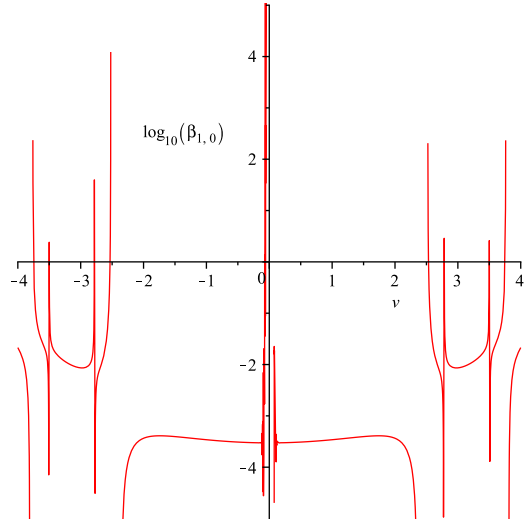}
\end{center} \caption{Behavior of the coefficient
$\beta_{1,0}$ in the method of PL$'$.}
\end{figure}
\begin{figure}
\begin{center}
\epsfxsize=8cm\epsfysize=8cm
\epsfbox{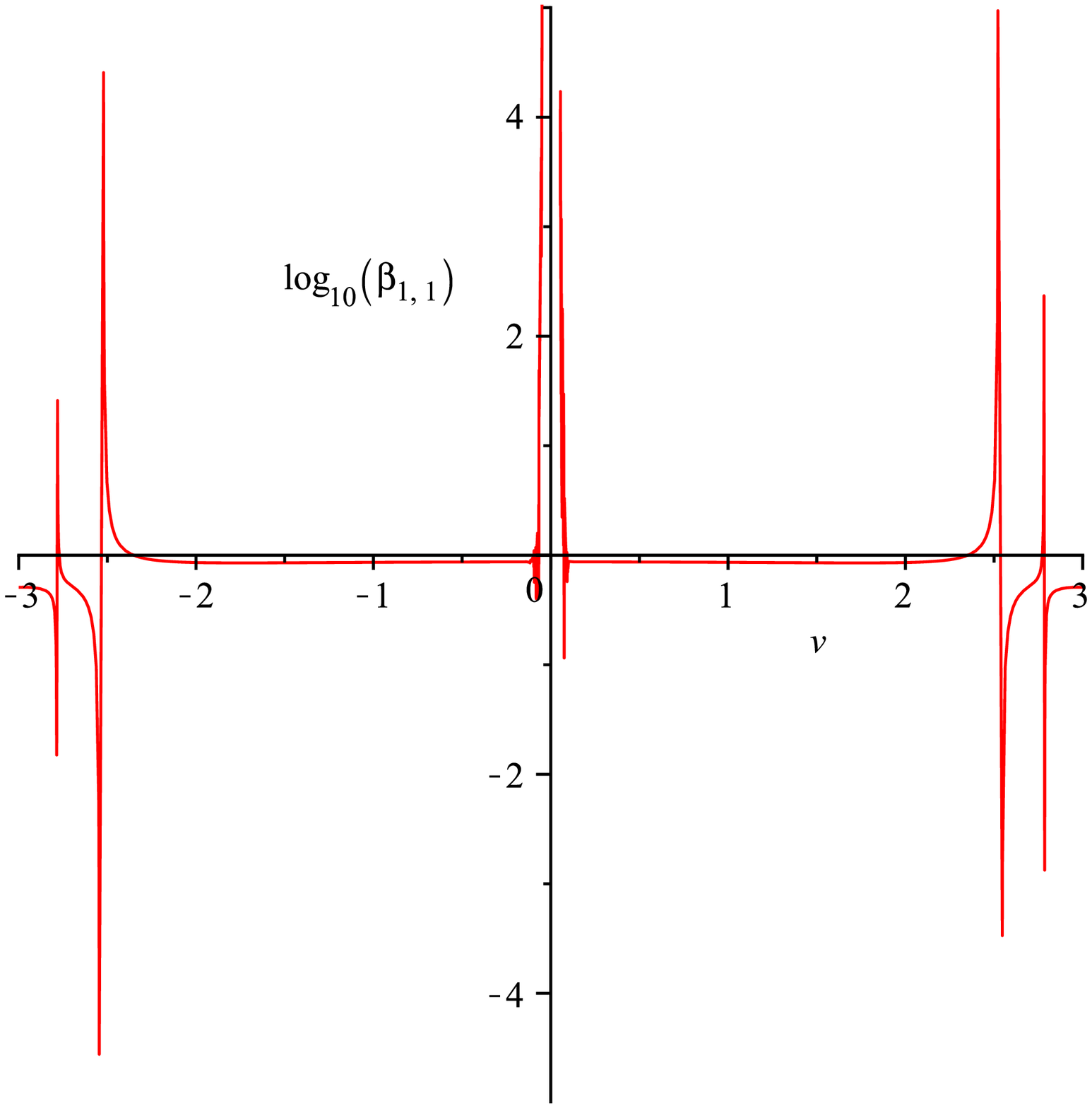}\end{center} \caption{Behavior of the coefficient
$\beta_{1,1}$ in the method of PL$'$.}
\end{figure}
\begin{figure}\begin{center}
\epsfxsize=8cm\epsfysize=8cm
\epsfbox{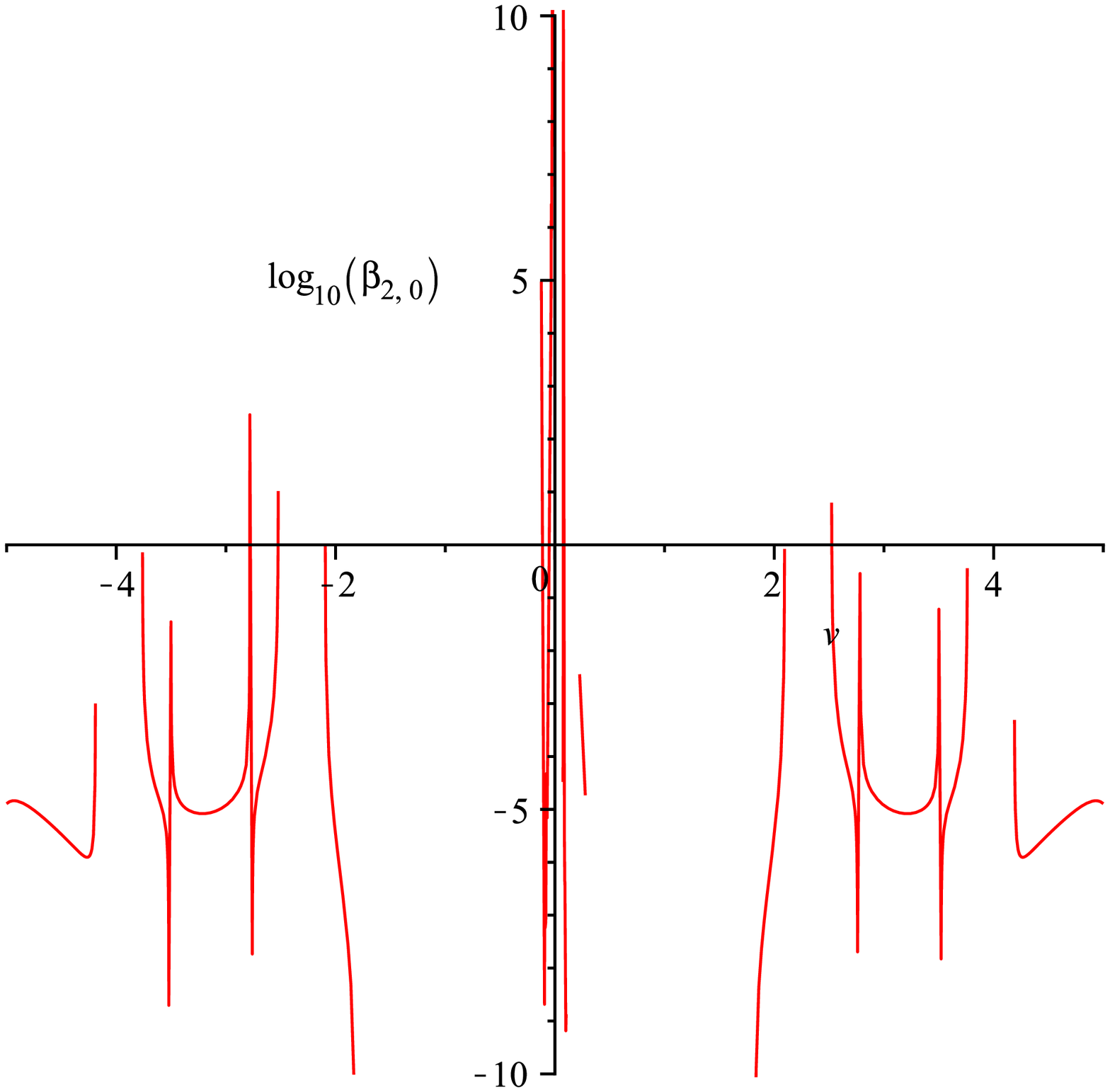}
\end{center} \caption{Behavior of the coefficient
$\beta_{2,0}$ in the method of PL$'$.}
\end{figure}
\begin{figure}\begin{center}
\epsfxsize=8cm\epsfysize=8cm
\epsfbox{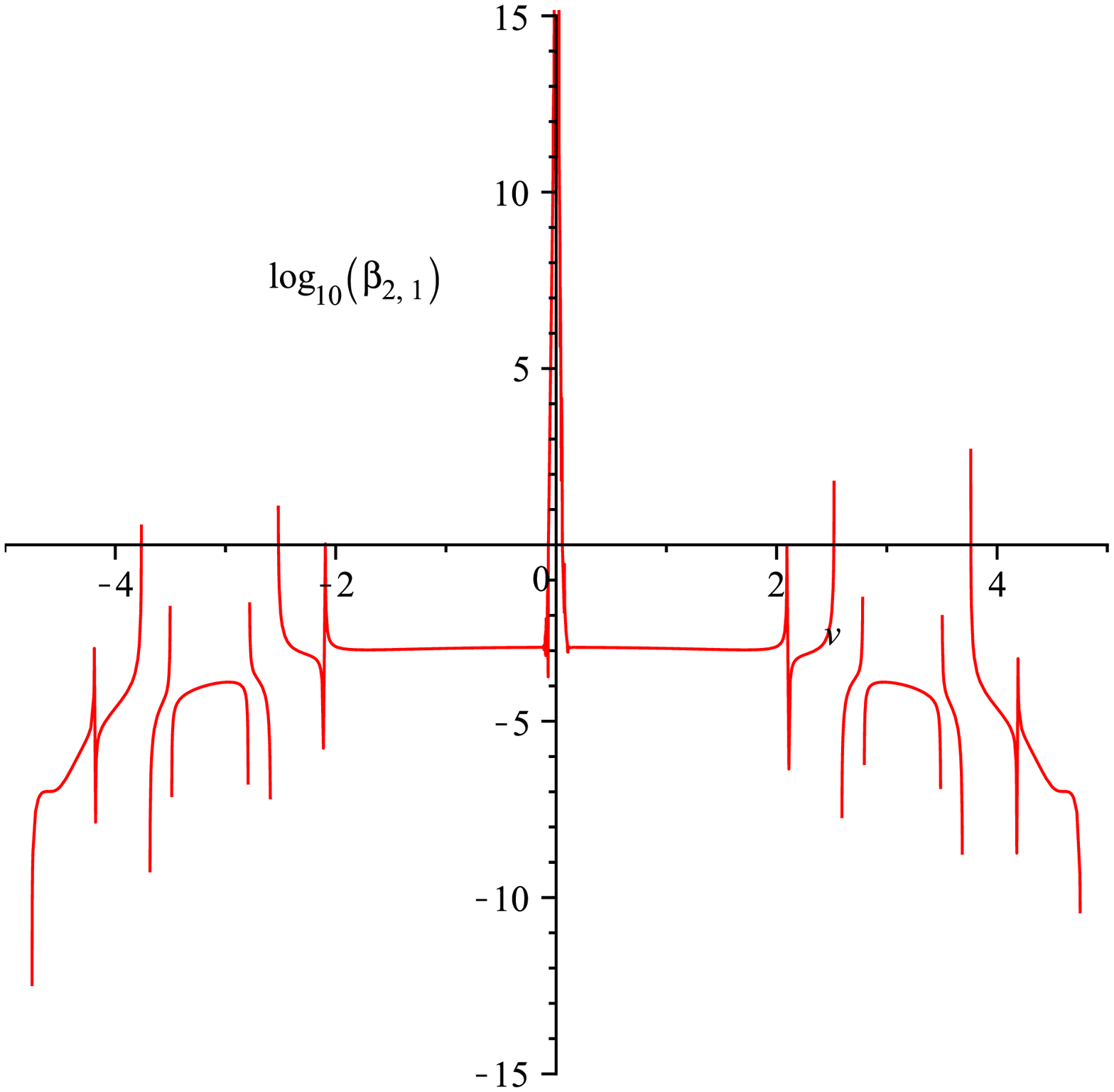}\end{center} \caption{Behavior of the coefficient
$\beta_{2,1}$ in the method of PL$'$.}
\end{figure}
\begin{figure}\begin{center}
\epsfxsize=8cm\epsfysize=8cm
\epsfbox{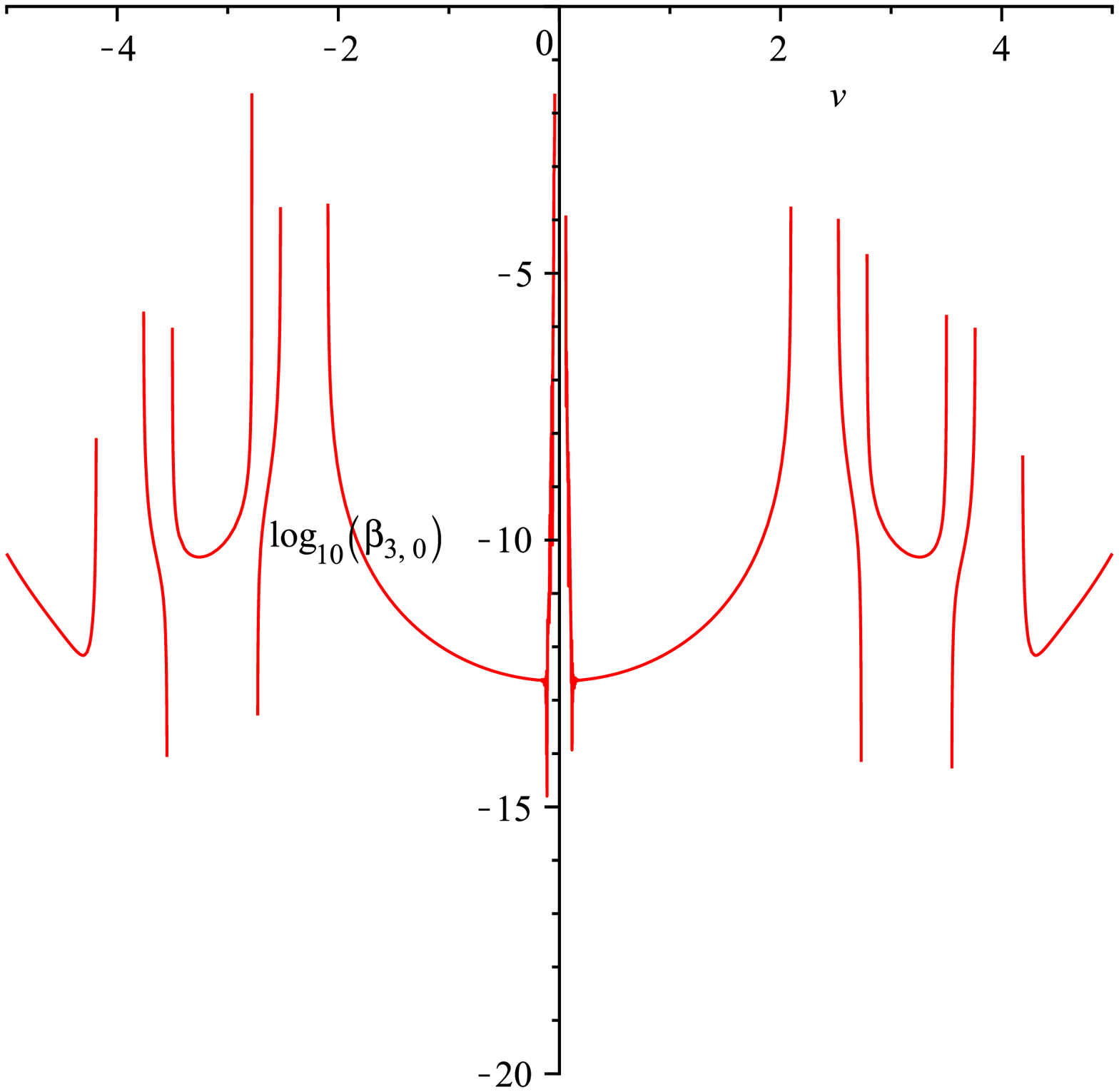}
\end{center} \caption{Behavior of the coefficient
$\beta_{3,0}$ in the method of PL$'$.}
\end{figure}
\begin{figure}\begin{center}
\epsfxsize=8cm\epsfysize=8cm
\epsfbox{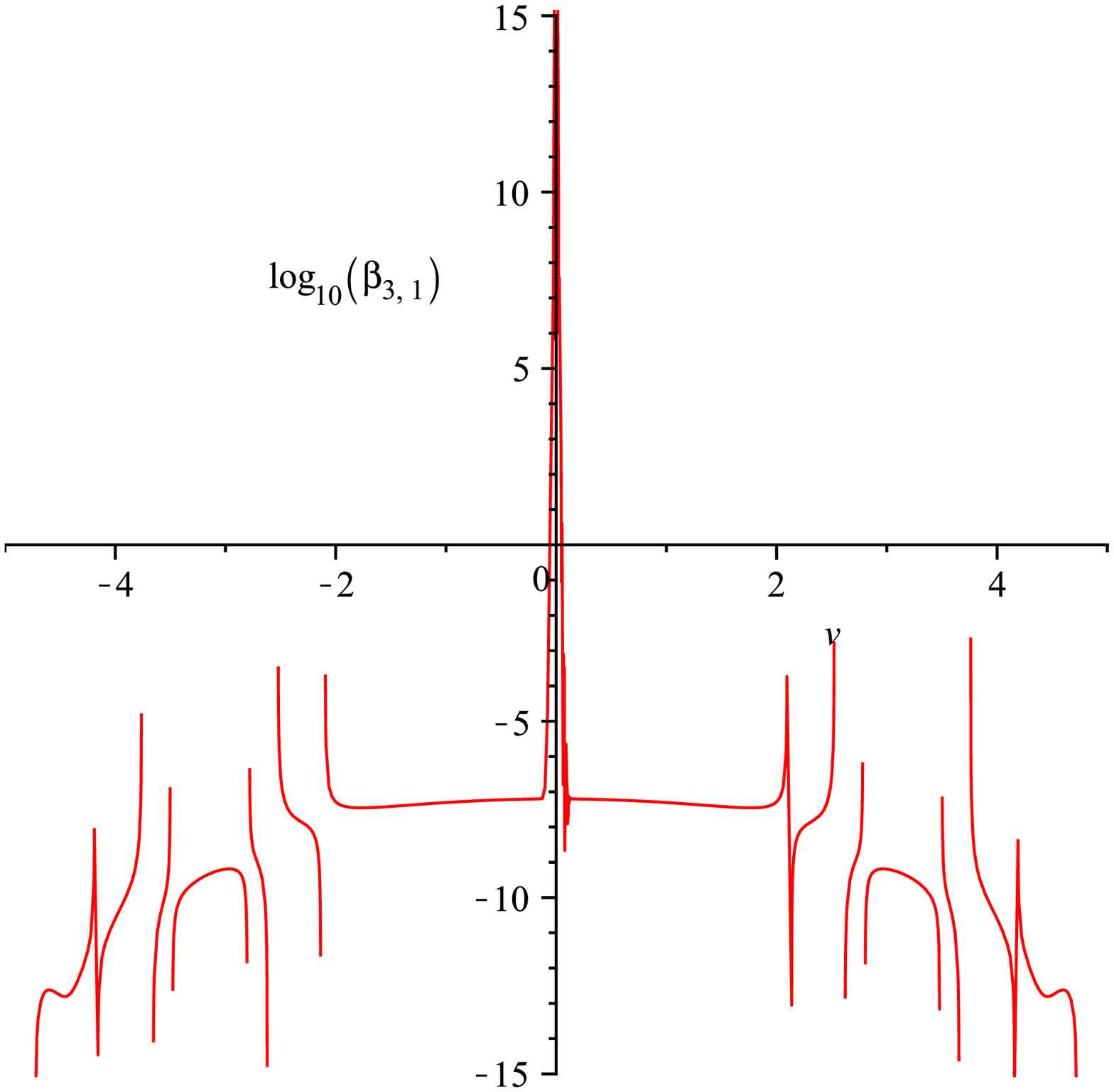}\end{center} \caption{Behavior of the coefficient
$\beta_{3,1}$ in the method of PL$'$.}
\end{figure}
%\\\\\\\\\\\\\\\\\\\\\\\\\\\\\\\\\\\\\\\\\\\\\\\\\\\\\\\\\\\\\\\\\\\\\\\\\\\\\\\\\\\\\\\\\\\\\\\\\\\\\\\\\\\\\\\\
\subsection{The second formula}
If $P=2$, $K=7$, so we called PL$''$, we have
\[\beta_{1,0}=\frac{89}{1878}-\frac{7560}{313}\beta_{3,1},\quad \beta_{1,1}=\frac{850}{939}+
\frac{15120}{313}\beta_{3,1},\quad\beta_{2,0}=-\frac{1907}{1577520}+\frac{330}{313}\beta_{3,1},\]
\[\beta_{2,1}=\frac{30257}{788760}+\frac{6900}{313}\beta_{3,1},\quad
\beta_{3,0}=\frac{59}{3155040}-\frac{13}{626}\beta_{3,1},\quad\beta_{3,1}=\frac{1}{1080}\frac{A}{B},\]
where
\begin{eqnarray}
A&=&-14400+213800\,\cos \left( 3\,v \right) {v} ^{4}\cos \left( v
\right) -36000\,\cos \left( 2\,v \right) \cos \left( v \right)
{v}^{2}+14400\,\cos \left( 3\,v \right) \cos \left(
v \right) \cos \left( 2\,v \right)\nonumber\\
&-&72000\,\cos \left( 3\,v \right) \cos \left( v \right)
{v}^{2}+20275\,\cos \left( 3\,v \right) {v}^{4} \cos \left( 2\,v
\right) -93600\,\cos \left( 3\,v \right) {v}^{2}\cos
\left( 2\,v \right)\nonumber\\
&+&9660\,\cos \left( 3\,v \right) {v}^{6}\cos \left( 2\,v \right)
+20832\,\cos \left( 3\,v \right) {v}^{6}\cos \left( v \right)
-10332\,\cos \left( v \right) {v}^{6}\cos \left( 2\,
v \right) +14400\,\cos \left( v \right)\nonumber\\
&-&14400\,\cos \left( 3\,v \right) \cos \left( 2\,v \right)
-14400\,\cos \left( 2\,v \right) \cos \left( v \right) +14400\,\cos
\left( 2\,v \right) +14400\,\cos \left( 3\,v
\right)\nonumber\\
&-&116475\,\cos \left( 2\,v \right) {v}^{4}\cos \left( v \right)
+100800\,\cos \left( 3\,v \right) \cos \left( v \right) \cos \left(
2\,v \right) {v}^{2}\nonumber\\
&+&29400\,\cos \left( 3\,v \right) \cos \left( v \right) \cos \left(
2\,v \right) {v}^{4}+720\, \cos \left( 3\,v \right) \cos \left( v
\right) \cos \left( 2\,v
\right) {v}^{6}+7200\,\cos \left( v \right) {v}^{2}\nonumber\\
&+&2875\,\cos \left( v \right) {v}^{4}+1830\,\cos \left( v \right)
{v}^{6}+28800\, \cos \left( 2\,v \right) {v}^{2}+99200\,\cos \left(
2\,v \right) {v}^{ 4}-46848\,\cos \left( 2\,v \right)
{v}^{6}\nonumber\\
&+&64800\,\cos \left( 3\,v \right) {v}^{2}-249075\,\cos \left( 3\,v
\right) {v}^{4}+88938\,\cos \left( 3\,v \right) {v}^{6}-14400\,\cos
\left( 3\,v \right) \cos \left( v \right)\nonumber\\
&-&810\,\cos \left( 3\,v \right) {v}^{8}\cos \left( 2 \,v \right)
+1296\,\cos \left( 3\,v \right) {v}^{8}\cos \left( v \right)
-486\,\cos \left( v \right) {v}^{8}\cos \left( 2\,v
\right),\nonumber
\end{eqnarray}
and
\begin{eqnarray}
B&=&240\,\cos \left( v \right) -81\,\cos \left( 2\,v \right) {v}^
{4}\cos \left( v \right) -240\,\cos \left( 3\,v \right) \cos \left(
v \right) -240\,\cos \left( 2\,v \right) \cos \left( v
\right)\nonumber\\
&+&96\, \cos \left( 3\,v \right) {v}^{4}\cos \left( v \right)
+75\,\cos \left( v \right) {v}^{4}-1107\,\cos \left( 2\,v \right)
\cos \left( v \right) {v}^{2}+240\,\cos \left( 3\,v \right) \cos
\left( v \right) \cos \left( 2\,v \right)\nonumber\\
&+&115\,\cos \left( v \right) {v}^{2}+992\,\cos \left( 3\,v \right)
\cos \left( v \right) {v}^{2}-240-15\,\cos \left( 3\,v \right)
{v}^{4}\cos \left( 2\,v \right) +115\,\cos \left( 3\,v \right)
{v}^{2}\cos \left( 2\,v \right)\nonumber\\
&-&240\,\cos \left( 3\,v \right) \cos \left( 2\,v \right) -480\,\cos
\left( 2\,v \right) {v}^{4}+992\,\cos \left( 2\,v \right)
{v}^{2}+405\,\cos \left( 3\,v \right) {v}^{4}-1107\,\cos \left( 3\,v
\right) {v}^{2}\nonumber\\
&+&240\,\cos \left( 3\,v \right) +240\,\cos \left( 2\,v
\right){v}^{6}.\nonumber
\end{eqnarray}
For small values of $v$ the above formulae are subject to heavy
cancelations. In this case the following Taylor series expansion
must be used:\\
\begin{eqnarray}
\beta_{1,0}&=&{\frac {229}{7788}}+{\frac
{318283}{657073560}}\,{v}^{2}+{\frac {
1512119}{118091281968}}\,{v}^{4}+{\frac {22946405723893}{
44038057651609507200}}\,{v}^{6}\nonumber\\
&+&{\frac {18296930817563773}{
651639946682396199939840}}\,{v}^{8}+{\frac {2913158423117216376847}{
1649365869047813021667729024000}}\,{v}^{10}\nonumber\\
&+&{\frac{8050460719799780764991137}{68936236116374773928415735195494400}}\,{v}^
{12} +\cdots,\nonumber
\end{eqnarray}
\begin{eqnarray}
\beta_{1,1}&=&{\frac {3665}{3894}}-{\frac
{318283}{328536780}}\,{v}^{2}-{\frac {
1512119}{59045640984}}\,{v}^{4}-{\frac {22946405723893}{
22019028825804753600}}\,{v}^{6}\nonumber\\
&-&{\frac {18296930817563773}{
325819973341198099969920}}\,{v}^{8}-{\frac {2913158423117216376847}{
824682934523906510833864512000}}\,{v}^{10}\nonumber\\
&-&{\frac{8050460719799780764991137}{34468118058187386964207867597747200}}\,{v}^
{12} +\cdots,\nonumber
\end{eqnarray}
\begin{eqnarray}
\beta_{2,0}&=&-{\frac {1}{2360}}-{\frac
{45469}{6021183168}}\,{v}^{2}-{\frac {
99714443}{586162181404800}}\,{v}^{4}-{\frac {4808531881}{
1130388645602810880}}\,{v}^{6}\nonumber\\
&-&{\frac {176305401655838711}{
1741655857496586207111936000}}\,{v}^{8}-{\frac
{76862259930526632407}{
35266441127276874790568169676800}}\,{v}^{10}\nonumber\\
&-&{\frac{1089463503416967799081153}{26321108335343095499940553438279680000}}\,{
v}^{12}+\cdots,\nonumber
\end{eqnarray}
\begin{eqnarray}
\beta_{2,1}&=&{\frac {711}{12980}}-{\frac
{1045787}{2365464816}}\,{v}^{2}-{\frac {
10184496007}{921111999350400}}\,{v}^{4}-{\frac {162691107254479}{
369919684273519860480}}\,{v}^{6}\nonumber\\
&-&{\frac {4589005587219802631}{
195491984004718859981952000}}\,{v}^{8}-{\frac
{582588392135442371849}{
395847808571475125200254965760}}\,{v}^{10}\nonumber\\
&-&{\frac{2446080156637919477841851381}{25176712320762960912986616332267520000}}
\,{v}^{12} +\cdots,\nonumber
\end{eqnarray}
\begin{eqnarray}
\beta_{3,0}&=&{\frac {127}{39251520}}+{\frac
{45469}{109175299200}}\,{v}^{2}+{\frac
{274576771}{7368895994803200}}\,{v}^{4}+{\frac {115636672827803}{
39837504460225215744000}}\,{v}^{6}\nonumber\\
&+&{\frac {76494288958873853}{
360908278162557895351296000}}\,{v}^{8}+{\frac
{455635060442806091167}{
30449831428575009630788843520000}}\,{v}^{10}\nonumber\\
&+&{\frac{136101812396019182508073199}{131285852101792282007800655206318080000}}
\,{v}^{12} +\cdots,\nonumber
\end{eqnarray}
\begin{eqnarray}
\beta_{3,1}&=&{\frac {127}{39251520}}+{\frac
{45469}{109175299200}}\,{v}^{2}+{\frac
{274576771}{7368895994803200}}\,{v}^{4}+{\frac {115636672827803}{
39837504460225215744000}}\,{v}^{6}\nonumber\\
&+&{\frac {76494288958873853}{
360908278162557895351296000}}\,{v}^{8}+{\frac
{455635060442806091167}{
30449831428575009630788843520000}}\,{v}^{10}\nonumber\\
&+&{\frac{136101812396019182508073199}{131285852101792282007800655206318080000}}
\,{v}^{12} +\cdots.\nonumber
\end{eqnarray}
The phase-lag and the local truncation error for the PL$''$ method
are given by
\[pl_{PL''}=-\frac{141797497314423651101}{7514399077966985427530263756800000}v^{20},\] and
\[LTE_{PL''}=-\frac {45469h^{14}}{1697361329664000}\left(49\omega^4y^{(10)}+
y^{(14)}+36\omega^6y^{(8)}+14\omega^2y^{(12)}\right),\] where
$v=\omega h$, $\omega$ is the frequency and $h$ is the step length.
The behavior of the coefficients of the PL$''$ method are shown in
Figures 4, 5, 6.
\begin{figure} \epsfxsize=8cm\epsfysize=8cm
\begin{center}\epsfbox{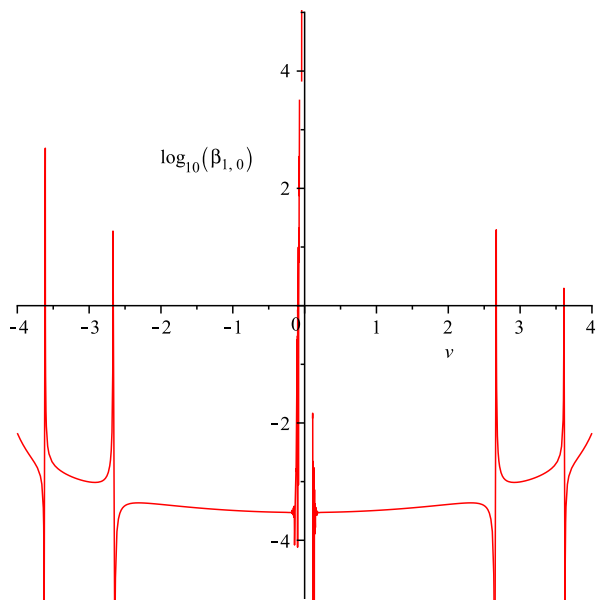}
\end{center} \caption{Behavior of the coefficient
$\beta_{1,0}$ in the method of PL$''$.}
\end{figure}
\begin{figure}
\begin{center}
\epsfxsize=8cm\epsfysize=8cm
\epsfbox{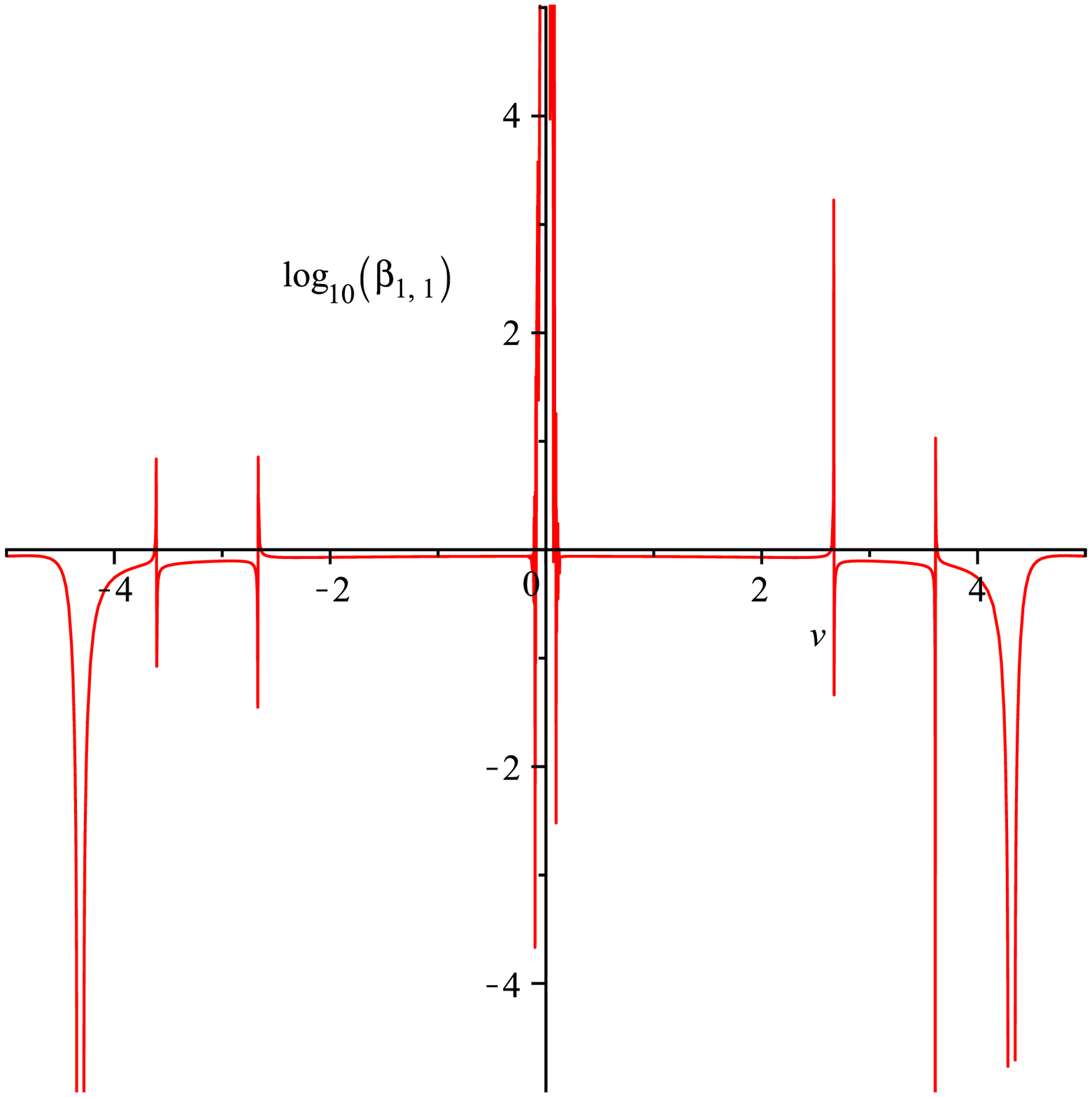}\end{center} \caption{Behavior of the coefficient
$\beta_{1,1}$ in the method of PL$''$.}
\end{figure}
\begin{figure}\begin{center}
\epsfxsize=8cm\epsfysize=8cm
\epsfbox{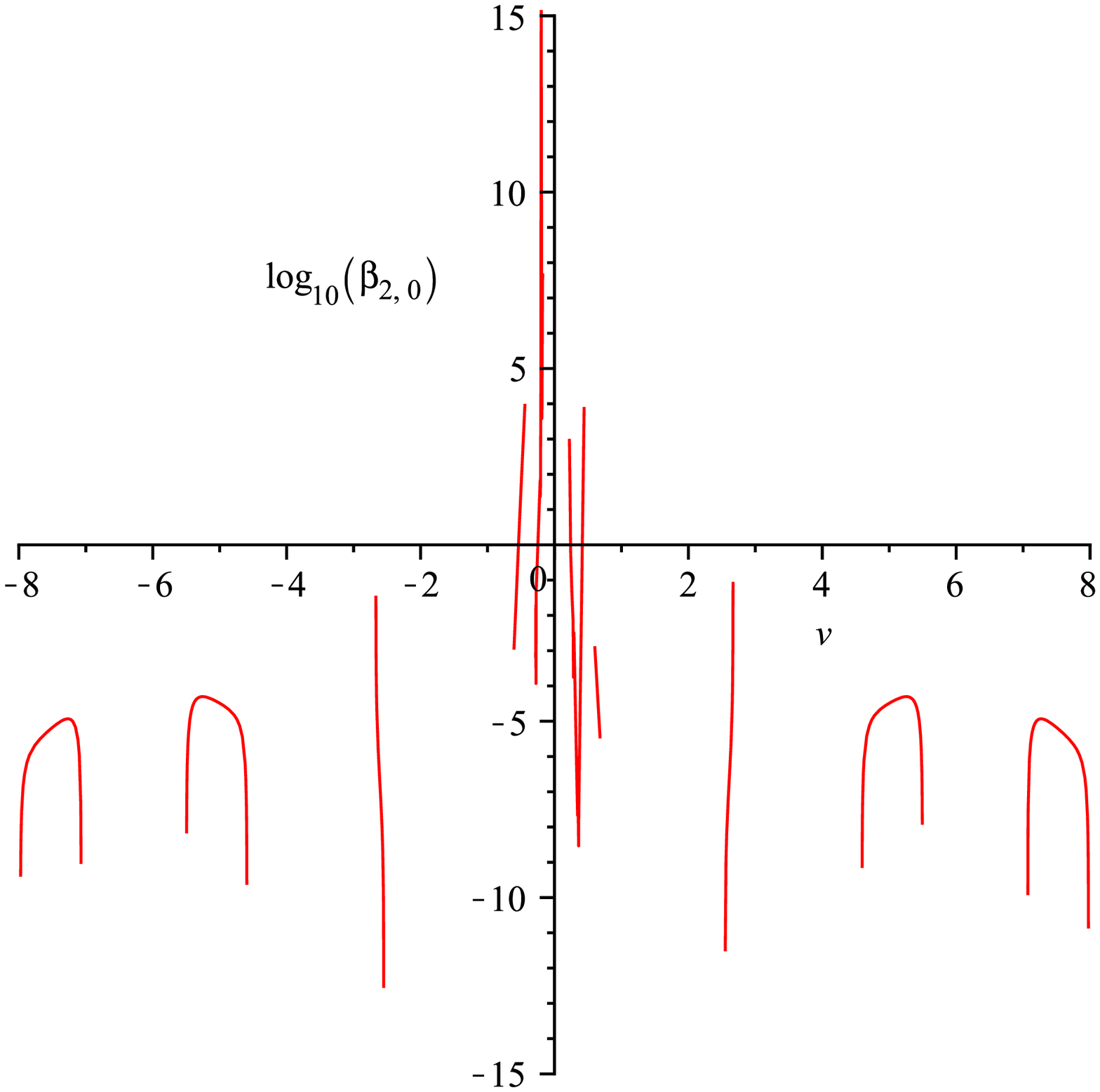}
\end{center} \caption{Behavior of the coefficient
$\beta_{2,0}$ in the method of PL$''$.}
\end{figure}
\begin{figure}\begin{center}
\epsfxsize=8cm\epsfysize=8cm
\epsfbox{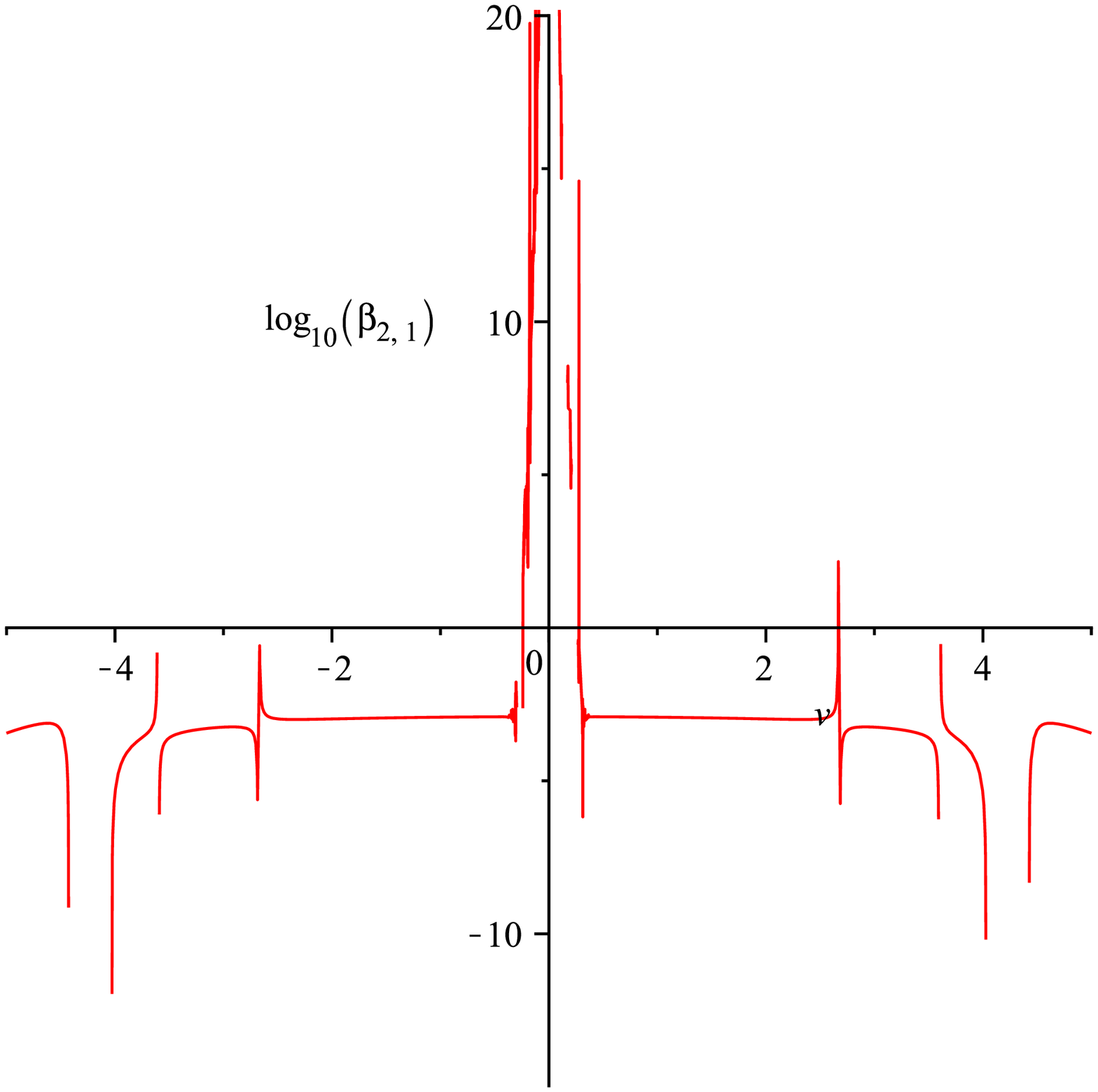}\end{center} \caption{Behavior of the coefficient
$\beta_{2,1}$ in the method of PL$''$.}
\end{figure}
\begin{figure}\begin{center}
\epsfxsize=8cm\epsfysize=8cm
\epsfbox{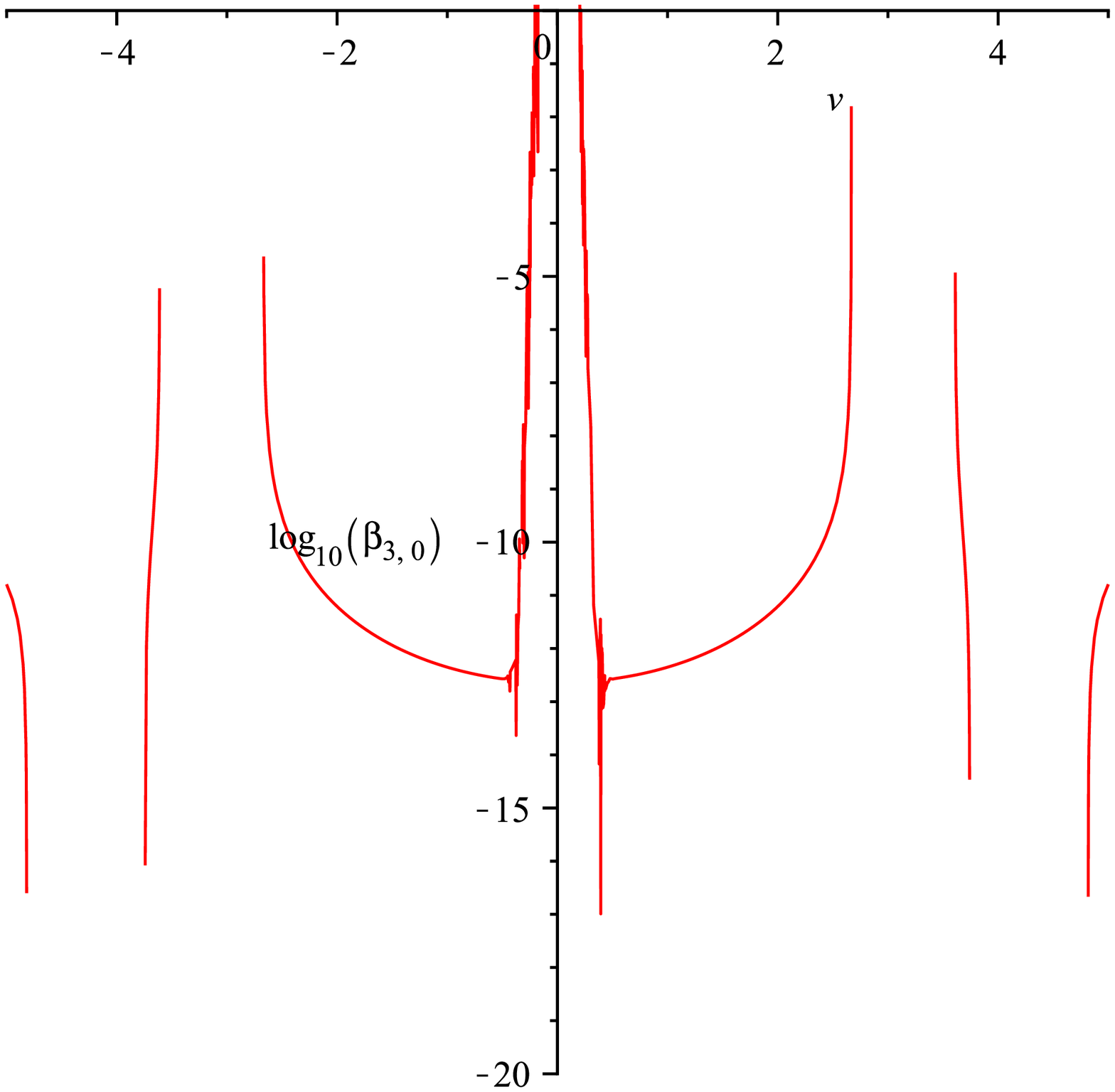}
\end{center} \caption{Behavior of the coefficient
$\beta_{3,0}$ in the method of PL$''$.}
\end{figure}
\begin{figure}\begin{center}
\epsfxsize=8cm\epsfysize=8cm
\epsfbox{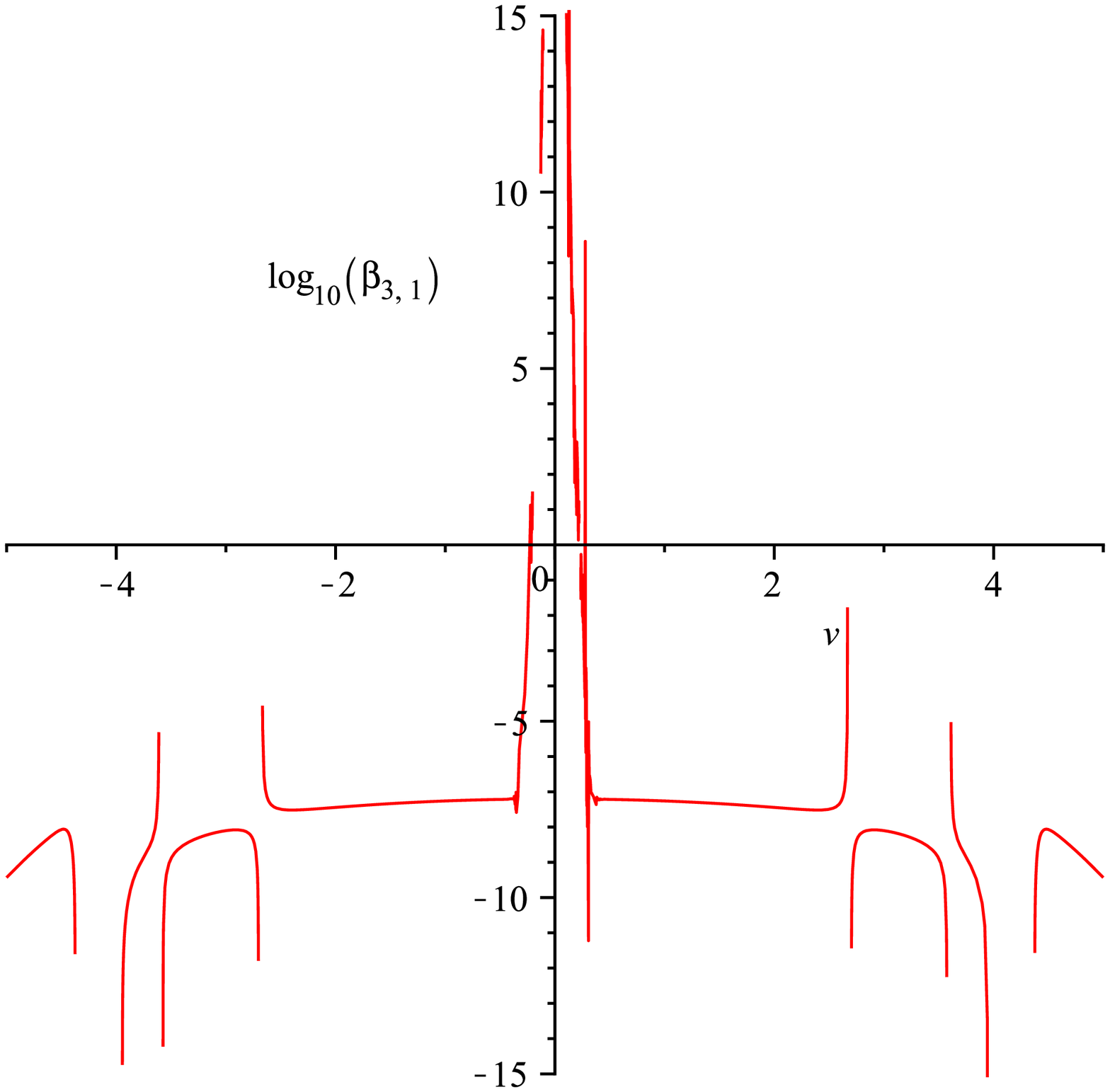}\end{center} \caption{Behavior of the coefficient
$\beta_{3,1}$ in the method of PL$''$.}
\end{figure}
The characteristic equation $\Omega(s:v^2)=A(v)s^2-2B(v)s+A(v)=0$
has complex roots of unit magnitude when
$\left|\cos(\theta(v))\right|=\left|\frac{B(v)}{A(v)}\right|<1$,
i.e. when $A(v)^2\pm B(v)^2>0$. Substituting for $A(v)$ and $B(v)$
for these the two-step methods, the interval of periodicity of the
classical Obrechkoff method, PL$'$ and PL$''$ methods when
$v\rightarrow0$ are obtained [0, 25.2004], [0,408.04] and
[0,1428.84] respectively.

%\\\\\\\\\\\\\\\\\\\\\\\\\\\\\\\\\\\\\\\\\\\\\\\\\\\\\\\\\\\\\\\\\\\\\\\\\\\\\\\\\\\\\\\\\\\\\\\\\\\\\\\\\\\\\\\\\\\\\\\\\\\\\\\\\\\\\\\\\\\
\section{Numerical example}
In this section, we present some numerical results obtained by our
new two-step trigonometrically-fitted Obrechkoff methods and compare
them with those from other multistep methods as\\
Achar: The 12th order Obrechkoff method of Achar \cite{Achar}.\\
Daele: The 12th order Obrechkoff method of Van Daele \cite{Da}.\\
Neta: The P-stable 8th-order super-implicit method of Neta \cite{Net}.\\
Simos: The 12th order Obrechkoff method of Simos \cite{Si}.\\
Wang: The 12th order Obrechkoff method of Wang \cite{Wan}.\\
%\\\\\\\\\\\\\\\\\\\\\\\\\\\\\\\\\\\\\\\\\\\\\\\\\\\\\\\\\\\\\\\\\\\\\\\\\\\\\\\\\\\\
\begin{example}
We consider the nonlinear \emph{undamped Duffing equation}
\begin{equation}\label{duf}
y''=-y-y^3+B\cos(\omega x),\quad y(0)=0.200426728067,\quad y'(0)=0,
\end{equation}
\end{example}
where $B=0.002$, $\omega=1.01$ and
$x\in\left[0,\frac{40.5\pi}{1.01}\right]$. We use the following
exact solution for \eqref{duf} from \cite{Ne},
\[g(x)=\sum_{i=0}^3K_{2i+1}\cos((2i+1)\omega x),\]
where
\begin{eqnarray}
\{K_1,K_3,K_5,K_7\}&=&\{0.200179477536,0.246946143\times10^{-3},\nonumber\\
&&0.304016\times10^{-6},0.374\times10^{-9}\}.\nonumber
\end{eqnarray}
In order to integrate this equation by a Obrechkoff method, one
needs the values of $y'$, which occur in calculating $y^{(4)}$.
These higher order derivatives can all be expressed in terms of
$y(x)$ and $y'(x)$ through \eqref{duf}, i.e.
\begin{eqnarray}
y^{(3)}(x)&=&-(1+3y^2(x))y'(x)-B\omega\sin(\omega x)\nonumber,\\
y^{(4)}(x)&=&-(1+3y^2(x))y''(x)-6y(x)y'(x)^2-B\omega^2\cos(\omega
x),\nonumber
\end{eqnarray}
The absolute errors at $x=\frac{40.5\pi}{1.01}$, for the new method,
in comparison with methods of Simos, Daele, Achar, Wang and the new
method are given in table 3.1 and the CPU times are listed in Table 3.2.
Also the absolute errors at $x=2\pi(4\pi)8\pi$, with $h=\frac{\pi}{12}$, for the new method PL$''$,
in comparison with methods Neta and the new
method are given in table 3.3.
\begin{table}
\[ \begin{tabular}{cccccc}
     &&   \\
    \hline   $h$           &        PL$''$       &    Simos             &     Daele           &      Achar     &      Wang \\\hline
    $\frac{M}{500}$        &    6.08953e-12      &   3.1486e-4          &   4.0560e-5         &     4.0919e-5  &     4.0831e-5\\
    $\frac{M}{1000}$       &    7.98859e-12      &   1.8069e-5          &   1.8733e-6         &     1.2708e-6  &     1.2678e-6\\
    $\frac{M}{2000}$       &    5.52149e-12      &   1.0752e-6          &   3.8355e-8         &     3.9420e-8  &     3.9327e-8\\
    $\frac{M}{3000}$       &    7.27826e-12      &   2.0873e-7          &   5.1344e-9         &     5.1801e-9  &     5.1678e-9\\
    $\frac{M}{4000}$       &    6.99211e-12      &   6.5463e-8          &   3.1876e-9         &     1.2324e-9  &     1.2308e-9\\
    $\frac{M}{5000}$       &    6.64542e-12      &   2.6673e-8          &   9.8900e-10        &     4.0911e-10 &     4.0741e-10\\
            \hline
  \end{tabular}\]
\caption{Comparison of the end-point absolute error in the
approximations obtained by using Methods: methods of Simos, Daele,
Achar, Wang and the new method for Example 3.1.}
\end{table}
\begin{table}
\[ \begin{tabular}{cccccc}
     &&   \\
    \hline   $h$           &    PL$''$     &    Simos             &     Daele       &      Achar     &      Wang  \\\hline
    $\frac{M}{500}$        &    1.453      &   1.437              &   1.484         &     1.188      &     1.406   \\
    $\frac{M}{1000}$       &    2.874      &   2.892              &   2.938         &     2.312      &     2.891   \\
    $\frac{M}{2000}$       &    6.267      &   6.233              &   6.36          &     4.812      &     6.236   \\
    $\frac{M}{3000}$       &    9.859      &   9.859              &   9.719         &     7.548      &     9.546   \\
    $\frac{M}{4000}$       &    13.424     &   13.548             &   13.39         &     9.986      &     13.063  \\
    $\frac{M}{5000}$       &    16.857     &   16.922             &   16.969        &     12.86      &     16.499  \\
            \hline
  \end{tabular}\]
\caption{CPU time for the example 3.1, are calculated for comparison
among four methods: methods of Simos, Daele, Achar, Wang and our new method PL$''$.}
\end{table}
\begin{table}
\[ \begin{tabular}{cccc}
     &&   \\
    \hline          $x$      &      CPU Time for PL$''$     &        PL$''$           &    Neta                   \\\hline
                    $2\pi$    &          0.03120020            &   6.06453e-14      &   2.53e-7          \\
                    $4\pi$    &          0.07800050            &   1.81249e-13      &   1.01e-6         \\
                    $6\pi$    &          0.09360060            &   3.45171e-13      &   2.25e-6         \\
                    $8\pi$    &          0.23400150            &   5.09481e-13      &   3.95e-6         \\
                    $10\pi$  &          0.28080180            &   6.24098e-13      &   6.05e-6          \\
                                  \hline
  \end{tabular}\]
\caption{Comparison of the end-point absolute error in the
approximations obtained by using Methods: Neta and the new method for Example 3.1.}
\end{table}
%\\\\\\\\\\\\\\\\\\\\\\\\\\\\\\\\\\\\\\\\\\\\\\\\\\\\\\\\\\\\\\\\\\\\\\\\\\\\\\\\\\\\
\begin{example}
Consider the initial value problem
\[y''=-100y+99\sin(x),\quad y(0)=1,\quad y^{\prime}(0)=11,\]
\end{example}
with the exact solution $y(t)=\sin(t)+\sin(10t)+\cos(10t)$.
This
equation has been solved numerically for $0\leq x\leq10\pi$ using
exact starting values. In the numerical experiment, we take the step
lengths $h=\pi/50$, $\pi/100$, $\pi/200$, $\pi/300$, $\pi/400$ and
$\pi/500$. In Table 3.4, we present the absolute errors at the
end-point and the CPU times are listed in Table 3.5.
\begin{table}
\[ \begin{tabular}{ccccc}
     &&   \\
    \hline   $h$               &           PL$''$         &    Simos              &     Daele            &      Achar      \\\hline
    $\frac{\pi}{50}$        &    1.76536e-26      &   3.0541e-11          &   1.2018e-11         &     5.7910e-13  \\
    $\frac{\pi}{100}$       &    4.50405e-30      &   2.2800e-13          &   7.3450e-13         &     5.7910e-13  \\
    $\frac{\pi}{200}$       &    1.90628e-34      &   4.3960e-13          &   8.6240e-13         &     1.3172e-12  \\
    $\frac{\pi}{300}$       &    4.60850e-37      &   2.1074e-12          &   2.6342e-12         &     1.9640e-12  \\
    $\frac{\pi}{400}$       &    6.28113e-39      &   1.3768e-12          &   2.9310e-12         &     4.7813e-12  \\
    $\frac{\pi}{500}$       &    2.23002e-40      &   6.4658e-12          &   2.8868e-12         &     7.5018e-12  \\
            \hline
  \end{tabular}\]
\caption{Comparison of the end-point absolute error in the
approximations obtained by using Methods: methods of Simos, Daele,
Achar and the new method for Example 3.2.}
\end{table}
\begin{table}
\[ \begin{tabular}{ccccc}
     &&   \\
    \hline   $h$                &       PL$''$          &    Simos             &     Daele            &      Achar      \\\hline
    $\frac{\pi}{50}$         &    0.2652017      &   0.1716011          &   0.2496016          &     0.187201    \\
    $\frac{\pi}{100}$       &    0.5772037      &   0.5148033          &   0.5304034          &     0.452403    \\
    $\frac{\pi}{200}$       &    1.1388073      &   0.8580055          &   0.8268053          &     0.748805    \\
    $\frac{\pi}{300}$       &    1.8096116      &   1.1388073          &   1.1544074          &     0.951606    \\
    $\frac{\pi}{400}$       &    2.496016       &   1.3884089          &   1.4040091          &     1.23241     \\
    $\frac{\pi}{500}$       &    2.9484189      &   1.7004109          &   1.7784114          &     1.46641     \\
            \hline
  \end{tabular}\]
\caption{CPU time for the example 3.2, are calculated for comparison
among four methods: methods of Simos, Daele, Achar and the new
method PL$''$.}
\end{table}
%\\\\\\\\\\\\\\\\\\\\\\\\\\\\\\\\\\\\\\\\\\\\\\\\\\\\\\\\\\\\\\\\\\\\\\\\\\\\\\\\\\\
\begin{example}
Consider the initial value problem
\[y''=\frac{8y^2}{1+2x},\quad y(0)=1,\quad y'(0)=-2,\quad x\in[0,4.5],\]
\end{example}
with the exact solution The theoretical solution of this problem is
\[y(x)=\frac{1}{1+2x}.\]
The absolute errors at $x=4.5$ for the new method, in comparison
with methods of Wang, Simos, Daele and Achar are given in the Table
3.6. The relative CPU times of computation of the new method in
comparison with the other four referred methods are given in Table
3.7.
\begin{table}
\[ \begin{tabular}{cccccc}
     &&   \\
    \hline   $h$            &   PL$''$        &    Simos              &     Daele            &      Achar      &  Wang      \\\hline
    $\frac{4.5}{500}$       &    2.74277e-21      &   1.2411e-7           &   1.2578e-7          &     1.2633e-7   &  1.2411e-7 \\
    $\frac{4.5}{1000}$      &    1.54818e-24      &   3.8166e-9           &   3.9035e-9          &     3.8481e-9   &  3.8166e-9 \\
    $\frac{4.5}{2000}$      &    5.84727e-28      &   1.1931e-10          &   1.2288e-10         &     1.2002e-10  &  1.1931e-10\\
    $\frac{4.5}{3000}$      &    5.22638e-30      &   1.9194e-11          &   2.0168e-11         &     1.4047e-11  &  1.9194e-11\\
    $\frac{4.5}{4000}$      &    1.78375e-31      &   7.8511e-12          &   7.8511e-12         &     2.6818e-12  &  7.8511e-12\\
    $\frac{4.5}{5000}$      &    1.28211e-32      &   1.6285e-12          &   1.6285e-12         &     7.4700e-14  &  1.6285e-12\\
            \hline
  \end{tabular}\]
\caption{Comparison of the end-point absolute error in the
approximations obtained by using five methods of Simos, Daele,
Achar, Wang and the new method for Example 3.3.}
\end{table}
\begin{table}
\[ \begin{tabular}{cccccc}
     &&   \\
    \hline   $h$            &   PL$''$      &    Simos         &     Daele        &      Achar    &  Wang      \\\hline
    $\frac{4.5}{500}$       &    0.3588023      &   0.359          &   0.343          &     0.187     &  0.312 \\
    $\frac{4.5}{1000}$      &    0.6084039      &   0.624          &   0.608          &     0.764     &  1.232\\
    $\frac{4.5}{3000}$      &    1.2792082      &   1.232          &   1.919          &     1.201     &  1.872\\
    $\frac{4.5}{4000}$      &    1.9344124       &   1.888          &   2.590          &     1.622     &  2.558\\
    $\frac{4.5}{5000}$      &    2.5584164         &   2.590          &   3.292          &     2.059     &  3.245\\
            \hline
  \end{tabular}\]
\caption{CPU time for the example 3.3, are calculated for comparison
among four methods of Simos, Daele, Achar, Wang and the new method PL$''$.}
\end{table}
%\\\\\\\\\\\\\\\\\\\\\\\\\\\\\\\\\\\\\\\\\\\\\\\\\\\\\\\\\\\\\\\\\\\\\\\
\section*{Conclusions}
In this paper, we have presented the new trigonometrically-fitted
two-step symmetric Obrechkoff methods of order 12. The details of
the procedure adapted for the applications have been given in
Section 2. With trigonometric fitting, we have improved the local
truncation error, phase-lag error, interval of periodicity and CPU
time for the classes of two-step Obrechkoff methods. The numerical
results obtained by the new method for some problems show its
superiority in efficiency, accuracy and stability.
%\\\\\\\\\\\\\\\\\\\\\\\\\\\\\\\\\\\\\\\\\\\\\\\\\\\\\\\\\\\\\\\\\\\\\\\\\\\\\\\\\\\\\\\\\\\\\\\\\\\\\\\\\\\\\\\
\section*{Acknowledgements}
The authors wish to thank the Professor Theodore E. Simos and the
anonymous referees for their careful reading of the manuscript and
their fruitful comments and suggestions.
%%%%%%%%%%%%%%%%%%%%%%%%%%%%%%%%%%%%%%%%%%%

\end{document}